\documentclass[10pt]{amsart}

\usepackage{verbatim}

\usepackage{amscd,amsmath,amssymb, hyperref}
\usepackage[mathscr]{eucal}
\usepackage{graphicx}

\theoremstyle{plain}

\numberwithin{equation}{section}

\newtheorem{theorem}{Theorem}[section]
\newtheorem{maintheorem}[theorem]{Main Theorem}
\newtheorem*{maintheorem*}{Main Theorem}

\newtheorem{corollary}[theorem]{Corollary}
\newtheorem{question}[theorem]{Question}

\newtheorem{definition}[theorem]{Definition}

\theoremstyle{remark}
\newtheorem{remark}{Remark}[section]
\newtheorem*{claim*}{Claim}
\newtheorem*{example*}{Example}

\newtheorem*{remark*}{Remark}
\newcommand{\R}{\mathbf{R}}

\DeclareMathOperator{\vol}{vol}

\begin{document}

\title[continuity of optimal transportation, positive curvature]
{Counterexamples to continuity of optimal transportation on
positively curved Riemannian manifolds}

\author{Young-Heon Kim
}

\address{Department of Mathematics, University of Toronto\\
  Toronto, Ontario Canada M5S 2E4}

\email{yhkim@math.toronto.edu
}
\date{\today}


\thanks{\copyright 2007 by the author. All rights reserved.}

\begin{abstract}
Counterexamples to continuity of optimal transportation on Riemannian manifolds with everywhere positive sectional curvature are provided. These examples show that the condition {\bf A3w} of Ma, Trudinger, \& Wang is not guaranteed  by positivity of sectional curvature.
\end{abstract}
\maketitle

\section{Introduction}
This paper addresses a question (see Question~\ref{Q:main}) in both optimal transportation theory and Riemannian geometry.
The question is explained in the following. For general notions we refer to the books by Villani \cite{V} \cite{V2} for optimal transport theory and the book by Cheeger and Ebin \cite{CE} for Riemannian geometry.

In optimal transportation, one considers two measure distributions $\rho$, $\bar \rho$ --- with the same total measure --- on domains $M$, $\bar M$, respectively, and one seeks for a minimizing (measurable) map $F:M \to \bar M$ for moving $\rho$ to $\bar \rho$ while it costs certain amount to move each unit mass at one location to another: this cost is given as a real valued function $c=c(x,\bar x)$ on the product $M \times \bar M$.

The case $c(\cdot, \cdot)=\mathop{\rm dist}^2(\cdot ,\cdot )/2$ for Riemannian distance $\mathop{\rm dist}$ on a Riemannian manifold $M = \bar M$ has been of great interest among researchers and an existence and uniqueness theory of optimal maps $F$ has been known for this case by the works of Brenier \cite{B} for Euclidean spaces and McCann \cite{M} for general Riemannian manifolds. Note that the distance squared cost $c$ (when differentiable) satisfies $\nabla_x c(x, \bar x) = (\exp_x )^{-1} (\bar x)$, and thus it can an be regarded as the canonical cost function for a Riemannian manifold --- when we say about a Riemannian manifold in this paper we always mean the manifold together with its distance squared cost.


The present work concerns the regularity of optimal transportation maps for Riemannian distance squared costs. A key notion is the so-called \emph{A3 weak condition} denoted as {\bf A3w} (see Definition~\ref{D:A3}). Ma, Trudinger, and Wang \cite{MTW}\cite{TW}\cite{TW2} have introduced and used this notion to develop a regularity theory of optimal transportation maps for general cost functions extending the results of Delano\"e \cite{D}, Caffarelli \cite{C}\cite{C2}, and Urbas \cite{U} for Euclidan distance squared costs. In fact, this {\bf A3w} is a necessary condition for continuity of optimal transport maps as shown later by Loeper \cite{L}:
he showed that if {\bf A3w} is violated then there exist smooth source and target measures $\rho$, $\bar \rho$
such that the optimal transportation map $F$ is not even continuous.
Moreover, for Riemannian distance squared costs, Loeper \cite{L} has shown that to satisfy {\bf A3w} the manifold should have nonnegative sectional curvature everywhere, and the standard round spheres $S^n$ satisfy {\bf A3w} (in fact a stronger condition so-called {\bf A3s}). This has led him to understand {\bf A3w} as a certain curvature condition, and to ask the following natural question.

\begin{question}\label{Q:main}{\bf (nonnegative curvature $\Longrightarrow$ A3w ?)}
Does every nonnegatively curved --- the sectional curvature is nonnegative everywhere --- Riemannian manifold satisfy {\bf A3w} for its distance squared cost?
\end{question}

As the main result of this paper, we answer Question~\ref{Q:main} \emph{negatively} by showing counterexamples.
\begin{maintheorem}\label{T:main}{\bf (nonnegative or positive curvature $\nRightarrow$ A3w)}
For each dimension $n$, there are complete (compact or noncompact) $n$-dimensional Riemannian manifolds with everywhere positive (nonnegative) curvature which do not satisfy  {\bf A3w}.
\end{maintheorem}
\begin{proof}
This result follows from Theorem~\ref{T:not A3w}, Corollary~\ref{C:positive curv}, and Remark~\ref{R:higher dim}: it is  shown that some shallow, smooth convex cones --- which are nonnegatively curved --- do not satisfy {\bf A3w}; by perturbation, positively curved examples are also obtained.
\end{proof}
As far as the present author knows, the examples we construct in this paper are the first  examples of nonnegatively or positively curved Riemannian manifolds where there are discontinuous optimal maps for smooth source and target measures $\rho$, $\bar \rho$. These examples confirm Trudinger's suspicion \cite{T} about Question~\ref{Q:main}.

\subsection*{Some perspectives on {\bf A3w} and Main Theorem~\ref{T:main}}
Let's first discuss  {\bf A3w} in some detail.  In its original form as introduced by Ma, Trudinger, \& Wang (to
show regularity of Monge-Amp\`ere type equations arising from optimal transportation theory), {\bf A3w} has been mysterious to researchers. The first geometric interpretation of {\bf A3w} is given by Loeper \cite{L} (see Theorem~\ref{T:local DASM}).
The present author and Robert McCann have given another more conceptual geometric interpretation \cite{KM2} by introducing a pseudo-Riemannian metric, say $h$, on the product space $M \times \bar M$ of the source and target domains. This metric $h$ is defined using the mixed second order partial derivatives of the cost function $c: M \times \bar M \to \R$ as the following non-degenerate\footnote{The non-degeneracy of $h$ needs the non-degeneracy of $D \bar D c$, a condition called {\bf A2} by Ma, Trudinger, \& Wang \cite{MTW} \cite{TW}, and each Riemannian distance squared cost on $M = \bar M$ (when differentiable) satisfies {\bf A2}.}  symmetric\footnote{$D \bar D c$ is the adjoint of $\bar D D c$.} bilinear form\footnote{This bilinear form is of type
$(+ + + \cdots - - - \cdots)$, i.e. it has the same number of positive and negative eigenvalues.} on $TM \oplus T\bar M$: \begin{equation}\label{metric}
h := \left(
\begin{array}{cc}
0 & -\frac{1}{2}\bar D D c \\
-\frac{1}{2}(D \bar D c) & 0
\end{array}
\right)
\end{equation}
where $D$, $\bar D$ denote the differentials of each $M$, $\bar M$, respectively. Then {\bf A3w} is equivalent to the following nonnegativity condition for the curvature of this pseudo-Riemannian metric $h$: namely,  for each $(x,\bar x) \in M\times \bar M$ and each tangent vector
$p \oplus \bar p \in  T_{(x, \bar x)}M =T_xM \oplus T_{\bar x} \bar M$ at $(x, \bar x)$,
\begin{align}\label{cross-curvature}
R_h((p\oplus 0) \wedge (0 \oplus \bar p), (p\oplus 0) \wedge (0 \oplus \bar p)) \ge 0
\hbox{ if
$h(p\oplus \bar p, p \oplus \bar p)=0$},
\end{align}
where $R_h$ denotes the curvature operator of $h$. The left-hand side of the inequality in \eqref{cross-curvature} is called \emph{cross-curvature}. Thus, {\bf A3w} can be interpreted as nonnegativity condition for cross-curvature of null-planes in $h$-geometry\footnote{A Riemannian manifold is said to be \emph{non-negatively cross-curved} if the inequality in \eqref{cross-curvature} holds without the condition $h(p\oplus \bar p, p \oplus \bar p)=0$.}. For the cost $c(x, \bar x) = \mathop{\rm dist}^2 (x, \bar x)/2$ on a Riemannain manifold $M = \bar M$,  $M$ is totally geodesically embedded as the diagonal of
$M \times M$ with respect to the pseudo-Riemannian metric $h$, and the cross-curvature in \eqref{cross-curvature} along this diagonal coincides with $R_M (p \wedge \bar p, p \wedge \bar p)$ where $R_M$ denotes the curvature operator of $M$
 (see \cite{KM2} for details).  This is another way to see Loeper's result \cite{L} that
{\bf A3w} implies nonnegative sectional curvature.

The pseudo-Riemannian metric $h$ and its curvature --- though they are local in the product space $M \times \bar M$ --- are global in nature with respect to the geometry of $M$ and $\bar M$: i.e. for a Riemannian manifold $M = \bar M$, local information concerning $h$ is equivalent to information about the global distance structure of $M$. Therefore {\bf A3w} --- the nonnegative cross-curvature condition for null-planes of $h$ --- is supposed to be a stricter restriction than nonnegative sectional curvature condition of $M$. Main Theorem~\ref{T:main} confirms this. As a consequence, this makes the following question of Trudinger \cite{T} much more interesting.
\begin{question}
For a Riemannian manifold $M$ with everywhere positive sectional curvature, let $R_M$ denote the curvature operator. Does there exist certain $\epsilon >0$ such that if
\begin{align*}
\|\nabla \log \| R_M \| \| \le \epsilon
\hbox{  \ with an appropriate point-wise norm $\| \cdot \|$},
\end{align*} then
$M$ satisfies {\bf A3w}?
\end{question}

Note that so far the only known examples of {\bf A3w} Riemannian manifolds\footnote{On the other hand, there are a lot of known examples of other types of cost functions satisfying {\bf A3w} \cite{MTW} \cite{TW}.} are modulo $C^\infty$-perturbations --- $C^2$ is maybe enough ---
the Euclidean space $\R^n$ (without perturbation), the standard $n$-dimensional sphere $S^n$ \cite{L}, and the Riemannian manifolds obtained from these by Riemannian coverings as considered by Cordero-Erausquin \cite{Co}, Delano\"e \cite{D2}, and Delano\"e and Ge \cite{DG}, more generally by Riemannian submersions\footnote{For example, the complex projective space $\mathbf{CP}^n$ with its Fubini-Study metric --- the sectional curvature $1 \le K_{\mathbf{CP}^n} \le 4$.} and products\footnote{For example, $M_1 \times \cdots \times M_k$, where $M_i=\R^l, S^m$, or $\mathbf{CP}^n$ for
$i=1, \cdots, k$.} as shown by the present author and McCann \cite{KM3}\footnote{In fact, we showed that (1) $S^n$ with its standard round metric is non-negatively cross-curved; (2) Riemannian submersions of {\bf A3w/ A3s} (resp. non-negatively cross-curved) Riemannian manifolds always induce  {\bf A3w /A3s} (resp. non-negatively cross-curved) Riemannian manifolds; (3) for products, if each factors are  non-negatively  cross-curved, then the resulting manifolds are non-negatively cross-curved, thus {\bf A3w} (but never {\bf A3s}); (4) and moreover, if one of the factors is not non-negatively cross-curved then the product is not {\bf A3w}. See \cite{KM2}\cite{KM3} for details and generalizations.}(see also \cite{KM2}). For all these unperturbed examples, $\| R_M \| = \mathop{\rm const}$.

\subsection*{Organizational remarks}
Although the notions and terminology in this paper have more general versions, they are
specialized to Riemannian distance squared costs for the sake of expositional simplicity.

This paper is organized as follows. In Section~\ref{S:pre} some preliminary notions and results are presented; Section~\ref{S:idea} explains the key idea of the counterexamples we construct; Section~\ref{S:Riemannian results} shows a Riemannian geometric result which is used in the main theorem; Section~\ref{S:example} is devoted to the construction of the counterexamples --- the main theorem of this paper.
\subsection*{Acknowledgment}
It is a great pleasure of the author to thank Robert McCann for his encouragement, remarks, and a lot of inspiring conversations, especially sharing his deep insights and knowledge.
He also thanks Philippe Delano\"e, Gregoire Loeper, Neil Trudinger, C\'edric Villani, and Xu-Jia Wang for their helpful discussions, comments, and recent preprints.  He is grateful to
Adrian Nachman for his generous support and interest. He thanks all of 2006-07 participants of
Fields Analysis Working Group, for the stimulating environment which they helped
to create. He thanks Jin-Whan Yim who introduced him to Riemannian geometry. Of course, he does not forget the great help from his wife, Dong-Soon Shim.

\section{Preliminaries}\label{S:pre}

In this section, some preliminary results, the definitions of {\bf A3w} and other key notions are presented.

First, let's recall a $2$-dimensional version of the famous theorem of Toponogov \footnote{The full version can be found in \cite{CE}.} which is essentially used in the proof of main theorem (Theorem~\ref{T:not A3w}).
\begin{theorem}\label{T:Toponogov}{\bf (Toponogov's comparison theorem)}
Let M be a complete $2$-dimensional Riemannian manifold with sectional curvature $K_M \ge H$, and let $M^H$ be the simply connected $2$-dimensional space of constant curvature $H$. Let $\gamma_i : [0,1] \to M$ and $\bar \gamma_i : [0,1] \to  M^H$, $i=1,2$, be minimal geodesic segments, i.e. they are unique geodesic segments connecting their end points.
Suppose that $\gamma_1 (0) = \gamma_2(0)$, $\bar \gamma_1 (0) = \bar \gamma_2 (0)$;
$\measuredangle(\dot \gamma_1 (0), \dot \gamma_2 (0))=\measuredangle(\dot{ \bar \gamma}_1 (0), \dot{\bar \gamma}_2 (0)) < \pi$,
where $\measuredangle$ denotes the angle between tangent vectors.
Assume $L[\gamma_i] = L[\bar \gamma_i]$, $i=1,2$, where $L$ denotes arc-length.
Then
\begin{align}\label{comparision}
\mathop{\rm dist}(\gamma_1(1), \gamma_2 (1)) \le \mathop{\rm dist}(\bar \gamma_1(1), \bar \gamma_2 (1)),
\end{align}
where $\mathop{\rm dist}$ denotes the Riemannian distance.
Moreover, if there exists a point $z$ on $\gamma_1 \cup \gamma_2 \subset M$ such that $K_M(z) > H$, then
the inequality \eqref{comparision} is strict.
\end{theorem}

The following key notions are specialized to Riemannian distance squared costs for the sake of expositional brevity: in fact, they have more general definitions \cite{MTW} \cite{TW} \cite{L} \cite{KM2}\cite{V2}.

\begin{definition}\label{D:A3}{\bf ($c$-segment, A3w, and local DASM)}
Let $M$ be a complete Riemannian manifold and let $c$ denote the Riemannian distance squared cost, i.e. $c(x,y) =\mathop{\rm dist}^2(x,y)/2$ for $x, y \in M$.
\begin{itemize}
\item {\bf ($c$-segment)}\cite{MTW} A curve $t \in [0,1] \to M$ is called a \emph{$c$-segment} with respect to $x$,
if $\bar x(t) = \exp_x (p + t \xi)$, for some $p, \xi \in T_x M$ and
$c(x, \bar x(t)) = | p + t\eta|^2$.

\item {\bf (A3w)}\cite{MTW}\cite{TW}
$M$ is said to satisfy {\bf A3w} if for any triple $(x, \bar x(t), \eta)$ of a point $x \in M$, a $c$-segment
$t \in [0,1] \to \bar x(t) = \exp{( p + t \xi)}$, $p, \xi \in T_x M$, and a tangent vector $\eta \in T_x M$ with $\eta \bot \xi$,
\begin{align}\label{c-curvature}
 \frac{d^2}{dt^2}\Big{|}_{t=0} [-D^2_{xx} c] (x, \bar x(t)) \ \eta \  \eta \ge 0,
\end{align}
where $D^2_{xx}$ denotes the Riemannian Hessian with respect to the first argument of $c$.
\begin{remark}\label{R:cost-curvature}
Loeper \cite{L} calls the left-hand-side of the inequality \eqref{c-curvature} \emph{cost-sectional curvature} --- he has shown it coincides with Riemannian sectional curvature when $x=\bar x(0)$.
\end{remark}
\item {\bf (local DASM)}\footnote{The name {\bf DASM} is an abbreviation of ``Double mountain above Sliding Mountain" \cite{KM1} \cite{KM2}.}\cite{L} (c.f. \cite{KM2})
$M$ is said to satisfy \emph{\bf local DASM} if for any $x \in M$ and any $c$-segment
$t \in [0,1] \to  \bar x (t)$ with respect to $x$, there exists a neighborhood $U$ of $x$ such that
the function $f_t (\cdot) = -c(\cdot, \bar x(t) ) + c(x, \bar x (t))$ satisfies
\begin{align}\label{DASM}
f_t (y) \le \max[ f_0 (y), f_1 (y)]
, 0 \le t \le 1, \ \ \forall y \in U.
\end{align}
\end{itemize}
\end{definition}

The notion {\bf local DASM} can be understood as a geometric interpretation of {\bf A3w} because of the following theorem which is originally due to Loeper \cite{L}.
\begin{theorem}\label{T:local DASM}{\bf (A3w $\Leftarrow \Rightarrow$ local DASM)}
Let $M$ be a complete Riemannian manifold.
$M$ satisfies {\bf A3w} if and only if $M$ satisfies  {\bf local DASM}.
\end{theorem}
\begin{proof}
($\Longrightarrow$) This direction can be easily verified by the elementary and geometric method in [\cite{KM1}, Section 6] --- this method is applied to more general cases \cite{KM2}; see also \cite{V2} for a modified proof.  See \cite{L} for an analytical proof using the main result of \cite{TW}.

($\Longleftarrow$) This direction is shown for more general case by Loeper \cite{L}  using Taylor expansion argument.
\end{proof}

\section{the key idea of counterexample}\label{S:idea}
In this section, we demonstrate our key idea of the counterexample which we shall construct in Theorem~\ref{T:not A3w}. We shall find such a situation that {\bf local DASM} is violated for the Riemannian distance squared cost of a nonnegatively curved Riemannian manifold --- by Theorem~\ref{T:local DASM}, {\bf A3w} then shall be violated, too.

Let $M$ be a Riemannian manifold and let $c$ denote the Riemannian distance squared cost function
$c(\cdot , \cdot) = \mathop{\rm dist}^2 (\cdot, \cdot)/2$. In the following discussion, we assume that $c$ is differentiable whenever necessary.
Let $t \in [0,1] \to \bar x(t)$, be a $c$-segment with respect to $x$. Thus there exist $p, \eta \in T_x M$
such that $\bar x(t) = \exp_x (p + t \eta)$ and $c(x, \bar x(t)) = |p+t \eta|^2$. Choose a tangent vector $\xi \in T_x M$ with $\xi \bot \eta$. This orthogonality shall be crucial.  Let $y= \exp_x (s_0 \xi)$ for a sufficiently small $s_0 >0$. Now suppose that there exist a point $\bar x(t_0)$ for $1/2 < t_0 <2/3$ and a sufficiently small open neighborhood $B$ of $\bar x(t_0)$ such that the Gaussian curvature $K$ satisfies $K \equiv 0$ on $M \setminus B$ and $K >0$ on $B$. ($M$ is nonnegatively curved.)
Further assume that the tangent vectors $p + t_0 \eta$ and $\xi$ are not collinear.
Then we see that
\begin{align}\label{euclidean at 0 and 1}
c(y, \bar x(0)) & = |s_0 \xi -p|^2 ,\\\nonumber
c(y, \bar x(0)) & = |s_0 \xi - p - \eta|^2.
\end{align}
Since $K >0$ near $\bar x(t_0)$ and $K \ge 0$ everywhere, by Toponogov's comparison theorem
(Theorem~\ref{T:Toponogov}),
\begin{align}\label{smaller than euclidean}
c(y, \bar x(t_0)) < |s_0 \xi - p -t_0 \eta|^2.
\end{align}
By the orthogonality  $\xi \bot \eta$,
the function
$\tilde f(t) = -| s_0 \xi - p -t\eta|^2 + | p + t \eta|^2$
is constant! Thus by \eqref{euclidean at 0 and 1},
\begin{align*}
f_0 (y) & = -c(y, \bar x(0)) + c(x, \bar x(0)) \\
&= -|s_0 \xi -p|^2 + |p|^2 \\
&= -|s_0 \xi - p - \eta|^2 + | p + \eta|^2\\
& = f_1 (y),
\end{align*}
and by \eqref{smaller than euclidean},
\begin{align*}
f_{t_0} (y) & > - |s_0 \xi - p - t_0 \eta|^2 + | p+ t_0 \eta|^2 \\
& = \tilde f(t_0) = \tilde f (0) = f_0 (y) = f_1(y) \\
&= \max[ f_0 (y), f_1 (y)].
\end{align*}
This violates {\bf local DASM}.

\section{Some results in Riemannian geometry}\label{S:Riemannian results}
In the following we prove some  technical results
(Theorem~\ref{T:injectivity radius} and Corollary~\ref{C:inj lower bound}) in Riemannian geometry; these results seem to be new and they are used in our construction of a nonnegatively curved manifold that does not satisfy {\bf A3w} (see Theorem~\ref{T:not A3w}).

First recall some definitions (c.f. \cite{CE}).
Let $M$ be a Riemannian manifold and $x$ be a point in $M$.
Let $\sigma$ be a geodesic from $x$, i.e. $\sigma = \exp_x (t \xi)$, $t \ge 0, \xi \in T_x$.
A point $y \in M$ is called a \emph{conjugate point} of $x$ along $\sigma$
if $y = \exp_x (t_0 \xi)$ and $\exp_x$ is singular at $t_0 \xi$.
A point $y \in M$ is called a \emph{cut point} of $x$
if either there are two distinct minimal geodesics from $x$ to $y$
or there is a unique minimal geodesic $\gamma$ from $x$ to $y$ and $y$
is a conjugate point of $x$ along $\gamma$.
The \emph{injectivity radius} $\mathop{\rm inj}_M (x)$
and \emph{conjugate radius} $\mathop{\rm conj}_M (x)$ at $x$ are defined as follows:
\begin{align*}
\hbox{\rm inj}_M (x) & = \inf \{ \mathop{\rm dist}(x, y) \ | \
\hbox{$y$ is a cut point of $x$ }\},\\
\hbox{\rm conj}_M (x) &= \inf \{ \mathop{\rm dist} (x, y) \ | \
\hbox{$y$ is a conjugate point of $x$}\},
\end{align*}
where $\mathop{\rm dist}(x,y)$ denotes the Riemannian distance between $x$ and $y$.
Note that $\mathop{\rm inj}_M (x) \le \mathop{\rm conj}_M (x)$ and if a geodesic $\sigma$
from $x$ to $z$ has length less than $\mathop{\rm inj}_M (x)$ then $\sigma$ is minimal.

The following result and its corollary are used later in Section~\ref{S:example}, but they have their own independent interests.
\begin{theorem}\label{T:injectivity radius}{\bf (injectivity radius $=$ conjugate radius)}
Let $M$ be a $2$-dimensional simply connected manifold and let $K$ denote its Gaussian curvature.
Suppose $\int_M K_+ d\vol  < \pi $ where $K_+$ denotes the positive part of $K = K_+ - K_-$.
Then for every $x \in M$, $\mathop{\rm inj}_M (x) =\mathop{\rm conj}_M (x)$.

\end{theorem}
\begin{proof}
Suppose $\mathop{\rm inj}_M (x) < \mathop{\rm conj}_M (x)$. It is easy to see that
\begin{itemize}
\item there is a point $y \in M$ such that $\mathop{\rm dist}(x, y) = \mathop{\rm inj}_M (x)$;
\item there are two distinct minimal geodesics, say $t \in [0,1] \to \gamma_i(t) = \exp_x ( t \xi_i)$,
$i=0,1$, from $x$ to $y$,
e.g. $\gamma_0 (1)=\gamma_1 (1)=y$, $\dot{\gamma}_0 (0) \ne \dot{\gamma}_1 (0)$;
\item $\exp_x$ is non-singular at $\xi_0, \xi_1 \in T_xM$.
\end{itemize}

First, the tangent vectors $-\dot{\gamma}_0 (1)$, $-\dot{\gamma}_1 (1)$ at $y$
have exactly the opposite direction, i.e. they form angle $\pi$.
If this is not the case, then there exists a tangent vector $\eta$ at $y$
which forms the same angle with $-\dot{\gamma}_0 (1)$, $-\dot{\gamma}_1 (1)$
and this angle is less than $\pi / 2$. Because of non-singularity of $\exp_x$ at $\xi_0$, $\xi_1$,
there exist two distinct continuous families of geodesics $\{\gamma_i^s \}$, $i=0,1$, for $0<s \ll 1$,
 such that each $\gamma_i^s$ connects $x$ to the  points $\exp_y (s \eta)$.
(In these families  $\{\gamma_i^s \}$ the initial velocity vectors of geodesics are close to
that of $\gamma_i$, $i=0,1$, respectively.) By the first variation formula of arc-length  (c.f. \cite{CE})
and the angle condition ($< \pi /2$), it is easy to see that for small $0 < s \ll 1$,
these geodesics have lengths smaller than the length of $\gamma_i$, $i=0,1$.
This contradicts that the length of $\gamma_i$, $i=0,1$, is the same as $\mathop{\rm inj}_M (x)$.

By simple connectedness of $M$, the geodesics $\gamma_i$, $i=0,1$,
joined together bound a domain $D$ which is a topological disk.
Suppose $D$ is oriented in such a way that $\partial D= \gamma_0 - \gamma_1$
(here the parametrization $\gamma_i (t) = \exp_x (t \xi _t)$, $i=0,1$,
give the orientations of $\gamma_0$ and $\gamma_1$.)
Let $\vartheta$ be the counter-clockwise angle from $-\dot{\gamma}_1(0)$
to $\dot{\gamma}_0 (0)$ at $x$.
By Gauss-Bonnet theorem and our assumption $\int_D K < \pi$,
\begin{align*}
2 \pi = 2\pi \chi (D)= \vartheta + \int_D K < \pi + \pi
\end{align*}
which is a contradiction! This completes the proof.
\end{proof}

\begin{corollary}\label{C:inj lower bound}{\bf (large injectivity radius)}
In addition to the assumptions in Theorem~\ref{T:injectivity radius}, further assume
that $K \le \delta$ everywhere for a fixed $\delta >0$. Then,
$\mathop{\rm inj}_M (x) \ge \frac{\pi}{\sqrt{\delta}}$ for every $x \in M$.
\end{corollary}
\begin{proof}
This follows from Theorem~\ref{T:injectivity radius}
and Rauch's comparison theorem (c.f. \cite{CE}).
\end{proof}

\section{Nonnegatively or positively curved Riemannian manifolds not satisfying {\bf A3w} condition}\label{S:example}
In the following, it shall be shown that a shallow, smooth convex cone (which is apparently nonnegatively curved) fits well into the situation in Section~\ref{S:idea} and it does not satisfy {\bf local DASM}, thus not {\bf A3w} by Theorem~\ref{T:local DASM}.

\begin{theorem}\label{T:not A3w}{\bf (nonnegative curvature $\nRightarrow$ A3w)}
A nonnegatively curved complete (open or closed) manifold does not necessarily satisfy  {\bf A3w}.
\end{theorem}

\begin{proof}
This theorem shall be proven by constructing two nonnegatively curved surfaces such that one is open, the other is closed, and both of them do not satisfy {\bf local DASM}.  These examples then do not satisfy {\bf A3w} condition (see Theorem~\ref{T:local DASM}).

Fix cartesian coordinates $(a,b)$ of $\R^2$ with the origin $O= (0,0)$.
Let $\theta(a, b)$ denote the polar angle of $(a,b)$ with respect to the origin which is counter-clockwise from the positive
$a$-axis. For example, $\theta(1,0)= 0$ and $\theta (0, 1) = \pi /2$.

Let $\vartheta$ be a sufficiently small positive number, i.e. $0< \vartheta \ll 1$, and define an infinite conical sector $\mathcal{C}_\vartheta$ by
\begin{align*}
\mathcal{C}_\vartheta =\R^2 \setminus \{ (a, b) \in \R^2 \ | \ \frac{3\pi}{2} - \vartheta <
\theta(a,b) \le \frac{3\pi}{2} + \vartheta \}.
\end{align*}
Let $B$ denote $B(O,1) \cap \mathcal{C}_\vartheta$, where $B(O, 1) $ is the open unit disk centered at $O$.
By identifying the two sides of $\partial \mathcal{C}_\vartheta$, we view this domain $\mathcal{C}_\vartheta$ as an infinite cone in $\R^3$ with conical angle $2\pi-2\vartheta$. It is easy to see that by only perturbing the metric inside $B$, this cone can be changed to a smooth surface $\Sigma_\vartheta \subset \R^3$,  in such a way that
\begin{itemize}
\item[(1)] the Riemannian metric of $\Sigma_\vartheta$ is radially symmetric with respect to the center (the point corresponding to $O$);
\item[(2)] the Gaussian curvature $K$ of $\Sigma_\vartheta$ as a function on $\mathcal{C}_\vartheta$ satisfies $K \equiv 0 $ on $\mathcal{C}_\vartheta \setminus B$;
\item[(3)] $0 < K < \frac{1}{10000}$ on $B$. Here, (3) is possible since $0< \vartheta \ll 1$.
\end{itemize}
For later use, it is important to note that by Theorem~\ref{T:injectivity radius} and Corollary~\ref{C:inj lower bound}, (2) \& (3) imply
\begin{itemize}
\item[(4)]
$\mathop{\rm inj}_{\Sigma_\vartheta} (z) > 314$ for every $z \in \Sigma_\vartheta$.
\end{itemize}

In the following, the cartesian coordinates of $\mathcal{C}_\vartheta$ shall be used to describe
 points in $\Sigma_\vartheta$. Let $c$ denote the Riemannian distance squared cost of the surface $\Sigma_\vartheta$. Let $x=(10, 10)$, $y=(10, 11)$ and let $t \in [0,1] \to \bar x(t)$ be the $c$-segment with respect to $x$ from  $\bar x (0) = (-10, 1/2)$ to $\bar x (1) = (10, 1/2)$, which is  just an exponential image (with respect to the metric of $\Sigma_\vartheta$) of a line segment
in the tangent space at $x$, i.e. $\bar x(t) = \exp_x (p + t\eta)$ for $p , \eta \in T_x \Sigma_\vartheta$; moreover,
$c(x, \bar x(t))= | p + t \eta|^2$, $0\le t \le 1$.  Note that our conditions (1), (2), (3), \& (4) ensure that $c$ is differentiable for any pair of points inside
$B(O, 100) \cap \mathcal{C}_\vartheta$, so for all relevant points in our consideration. These conditions also make it clear that
\begin{itemize}
\item each point $\bar x(t)$ is connected to $x$ and $y$ by unique minimal geodesics;
\item the unique minimal geodesics from $x$ and $y$ to $\bar x(0)$, $\bar x(1)$ are the straight line segments outside $B$;
\item  the curve $t \in [0,1] \to \bar x(t)$ coincides with the straight line segment from $\bar x(0)$ to $\bar x(1)$ until it hits the ball $B$.
\end{itemize}
Therefore, there exists $0< t_0 <1$ with $\bar x(t_0)$ in $B$. Thus, by following the same lines of Section~\ref{S:idea}, $\Sigma_\vartheta$ does not satisfy {\bf local DASM}. This $\Sigma_\vartheta$ furnishes an example of open nonnegatively curved manifold not satisfying {\bf A3w}.

To get a closed surface example, first cut off a large geodesic ball $B_1$ of $O$, e.g. with radius $10000$, from the surface $\Sigma_\vartheta$, then glue a flat disk to $B_1$ along $\partial B_1$ and round-off the curve where the disk and $B_1$ are glued, in such a way that the resulting surface is smooth, radially symmetric from $O$, and convex (thus, the Gaussian curvature is nonnegative).
This completes the proof.
\end{proof}

\begin{corollary}\label{C:positive curv} {\bf (positive curvature $\nRightarrow$ {\bf A3w})}
A positively curved complete (open or closed) manifold does not necessarily satisfy {\bf A3w}.
\end{corollary}
\begin{proof}
It is possible to perturb (radially symmetrically) the above nonnegatively curved examples so that the resulting manifolds have
positive curvature everywhere. If {\bf local DASM} is violated at some points in the original manifolds, then it should be violated in the perturbed manifolds as well for sufficiently small perturbations. By Theorem~\ref{T:local DASM}, the corollary follows.
\end{proof}

\begin{remark}\label{R:higher dim}{\bf (higher dimensional examples)}
The examples in Theorem~\ref{T:not A3w}  are radially symmetric and we can easily construct higher dimensional radially symmetric examples in which our $2$-dimensional examples are isometrically and totally geodesically embedded. Then the higher dimensional examples do not satisfy {\bf A3w}, neither their positively curved radially symmetric perturbations. One may also consider taking Riemannian product of a positively curved but non-{\bf A3w} manifold with other positively curved  manifold, then certainly the resulting manifold is nonnegatively curved but violates {\bf A3w}; however,
it is not clear whether we can perturb the product to a positively curved manifold --- it is a famous conjecture of H. Hopf that $S^2 \times S^2$ does not carry a positive curvature metric. See \cite{KM2} \cite{KM3} for more consideration on products.
\end{remark}



\bibliographystyle{plain}

\end{document}